# The Lotka-Volterra Predator-Prey Model with Disturbance

**Ogethakpo A. J. and Ojobor S. A.**

Department of Mathematics, Delta State University Abraka, Delta State, Nigeria.
Email: josephogethakpo@gmail.com

The forces which drive growth, development, survival and change within an ecological system involving a predator and prey specie are not easily addressed in the field. To better understand the dynamics in the system, ecologists have turned to mathematical models. The predator-prey dynamics of rat and cat population in a given ecology is studied. The mathematical model proposed by Alfred J. Lotka and Vito Volterra called the Lotka-Volterra model for studying predator-prey dynamics is utilized. Assumptions were made to suit the given ecology. These assumptions lead to the modification of the Lotka-Volterra equations. The equilibrium and stability properties of the modified model is established. Results were simulated using MATLAB.

## INTRODUCTION

The interactions between individuals in a community can take distinct form. An interaction between two individuals may be beneficial to both individuals or the interaction may benefit one individual to the detriment of the other individual. When the interaction between two individual organisms benefits one to the detriment of the other, it is called an antagonistic interaction. An example is predation.

Predation is a biological interaction between two organisms, the predator and the prey. It is the act of the predator feeding on the prey, which may lead to the death of the prey. Predation provides energy to prolong the life and promote the reproduction of the predator to the detriment of the prey, which can have a major effect on the density and size of a population of organism as applied to population that when the death rate exceeds the birth rate in a population, the size of the population usually decreases.

Predation influences organisms at two ecological levels. At the individual level, the prey organism has an abrupt decline in fitness as measured by its lifetime reproduction success. This is because it will certainly not reproduce again. And at the community level, predation reduces the number of individuals in the prey population. But a decrease in the prey's population in turn affects the predator's population.

Predation can be a powerful determinant of a community's structure. It has the capacity to dynamically influence the numbers and quality of both the predator and the prey, as it acts as an important agent of natural selection on both groups. The best-known examples of predation involve carnivorous interactions in which one animal consumes another. e.g lion hunting antelope, wolves hunting moose, owls hunting mice, or shrews hunting worms and insects. Predation can also occur as parasitism, cannibalism, herbivory.

From the above exposition, one can say that no simple relationship exists between a community and its individual species. As was noted by Berryman (1992), the dynamic relationship between predators and their prey has long been and will continue to be one of the dominant themes in both ecology and mathematical ecology due to its universal existence and importance and depending on their specific settings of applications, predators and preys can take the forms of resource-consumer, plant-herbivore, parasite-host, tumor cells (virus)-immune system, susceptible-infectious interactions, etc. Since the complex dynamics for interactive species are common in the real world, many researchers have investigated the processes that affect the dynamics of predator-prey relation and wanted to know what models can best represent species interactions.

The Lotka-Volterra model is the simplest model of predator-prey interactions. It is a pair of differential equations describing the interaction between a predator and a prey. Since the development of the first Predator-Prey Model called the Lotka-Volterra model, many researchers have worked on extending and modifying the model to include other factors which could also come to play in the course of the interaction between a predator and a prey. Aside from the different variations of the first Predator-Prey Model, various researchers have also attempted to apply the predator-prey equations and its numerous variations to different systems either for the purpose of estimating the quantitative behavior of the system or comparing the quantitative results from the model to observational data. The various applications of the Predator-Prey Model ranges from biological systems, ecological systems, fishery, economics, mobile operating systems, etc.

Based on available literatures, it is observed that the effect of disturbances in the ecological system of predator-prey interaction has not been modeled and analyzed. This study model this factor in an ecosystem in which two species (cat and rat) interact. Assumptions/conditions for existence of both species shows that, since human beings predate on the rats, it will affect the mean number of the prey and predator because that action reduces the quantity of food meant for the predator, this could make them starve to death. The concept of the generalized Predator-Prey Models that has been stated is extended for the model. The model like other Predator-Prey Models deals with general loss-win interactions and hence may have applications outside of the ecosystem.

## METHODS
### The basic idea of the model
Here, we briefly recall the basic idea of the Lotka-Volterra model which will be the guide for the following sections;

i) The predator species is totally dependent on the prey species as its only food supply.
ii) The prey species has an unlimited food supply.
iii) There is no threat to the prey's growth other than the specific predator.
iv) The rate at which the predator encounters the prey is jointly proportional to the sizes of the two populations.
v) A fixed proportion of encounter leads to the death of preys.

These assumptions lead to the Lotka-Volterra Predator-Prey Model:

$$\frac{dR}{dt} = a_1 R(t) - a_2 R(t)C(t)$$

$$\frac{dC}{dt} = -b_1 C(t) + b_2 R(t)C(t)$$

Where $a_1$ = natural growth rate of the rat, $a_2$ = carrying capacity of the rat, $b_1$ = natural growth rate of the cat, $b_2$ = carrying capacity of the cat.

Given the certain dynamics of the system under study, we also have the following assumptions:

i) C(r) denotes the predator response function. The models C(r) of functional response are assumed to be continuously differentiable on [0, ∞] and satisfy $C(0) = 0$, $C^1(r) > 0$ and $\lim_{r \to \infty} C(r) = k < \infty$

Note: A functional response in ecology is the intake rate of a consumer as a function of food density. It is associated with the numerical response, which is the reproduction rate of a consumer as a function of food density. Following C. S. Holling, functional responses are generally classified into three types, which are called Holling's type I, II, and III.

i) The Holling type-III function satisfies the assumptions made and will be used for the model.

ii) We will assume that the primary loss of rats is due to predation by the cats and human beings.

Mathematically, this is given by a negative term $-a_2 R(t)C(t)$ and $-a_3R(t)$. Here, $a_3$ is the total number of prey, n minus the number predated by human being, and m, that is, (n–m).

Combining these terms, we have the growth model for the rat population:

$$\frac{dR}{dt} = a_1 R(t) - a_2 R(t)C(t) - (n-m)R(t)$$

$$\frac{dR}{dt} = a_1 R(t) - a_2 R(t)C(t) - a_3 R(t) \quad \text{- - - - - -}$$

This is the prey's Equation (a)

For the predator, the growth of the Cat population can be expressed as $b_2 R(t)C(t)$. The human predation is negative on the Cat because it reduces the quantity of food meant for the predator, that is, $(n-m)C(t) = b_3C(t)$, where $b_3$ represents the total number of rats, n minus the number predated by human beings, m. The loss of Cat is presumed to be a type of reverse growth. Thus, in the absence of rats, the cat population declines in population to their own population which is expressed by the negative modeling term as $-b_1C(t)$.

The growth model for the Cat population gives:

$$\frac{dC}{dt} = -b_1 C(t) + b_2 R(t)C(t) - (n-m)C(t)$$

$$\frac{dC}{dt} = -b_1 C(t) + b_2 R(t)C(t) - b_3 C(t) \quad \text{- - - - - -}$$

This is the predator's Equation (b)

The model ignores the role of climate variation and the interactions of other species. Other significant factors ignored are the ages of the animals and the spatial distribution. The two differential equations above are intertwined into a system of differential equations with each growth model depending on the unknown variable (population) of the other.

**The model**
From the predator and prey's equation above, let

$a_1 = ur$ (natural growth rate of the rat)

$a_2 = \dfrac{r}{m}$ (carrying capacity of the rat)

$a_3 = b_3 = qEr$ (the predating coefficient/total effort on predating on the rat population)

$b_1 = vc$ (natural growth rate of the cat)

$b_2 = \dfrac{c}{n}$ (carrying capacity of the cat)

Thus, the rate of change of rats population with respect to time,

$$\frac{dr}{dt} = ur\left(1 - \frac{r}{m}\right) - \frac{kr^2 c}{p + r^2} - qEr \quad (i)$$

Similarly, the rate of change of cats population with respect to time,

$$\frac{dc}{dt} = vc\left(1 - \frac{c}{n}\right) + \frac{ekr^2 c}{p + r^2} - dc \quad (ii)$$

Where $\dfrac{kr^2}{p + r^2}$ is the predator's functional response k and p are positive constants with k being the maximum growth rate of the species and p the saturation constant; e is the conversion rate; d the death rate of the cat; r and c stands for rat and cat respectively.

In summary, incorporating disturbances and the Holling type III functional response, the mathematical model of the predator-prey system becomes:

$$\frac{dr}{dt} = ur\left(1 - \frac{r}{m}\right) - \frac{kr^2 c}{p + r^2} - qEr \quad (1)$$

$$\frac{dc}{dt} = vc\left(1 - \frac{c}{n}\right) + \frac{ekr^2 c}{p + r^2} - dc \quad (2)$$

We write the Equations 1 and 2 in their non-dimensionalized form in other to reduce the number of parameters. From Equation (1) when we set

$$\bar{r} = \frac{r}{\sqrt{p}} \Rightarrow r = \bar{r}\sqrt{p}, \qquad dr = \sqrt{p}\,d\bar{r},$$

we have

$$\frac{\sqrt{p}\,d\bar{r}}{dt} = u\bar{r}\sqrt{p}\left(1 - \frac{\bar{r}\sqrt{p}}{m}\right) - \frac{k\bar{r}^2 c\sqrt{p^2}}{p + \bar{r}^2\sqrt{p^2}} - qE\bar{r}\sqrt{p} \quad (3)$$

setting

$$\bar{t} = \frac{ut\sqrt{p}}{m} \Rightarrow t = \frac{m\bar{t}}{u\sqrt{p}}, \qquad dt = \frac{m\,d\bar{t}}{u\sqrt{p}},$$

we have

$$\frac{\sqrt{p}\,d\bar{r}}{\frac{m\,d\bar{t}}{u\sqrt{p}}} = u\bar{r}\sqrt{p}\left(1 - \frac{\bar{r}\sqrt{p}}{m}\right) - \frac{kpc\bar{r}^2}{p + p\bar{r}^2} - qE\bar{r}\sqrt{p}, \quad (4)$$

$$\frac{u\sqrt{p}^2\,d\bar{r}}{m\,d\bar{t}} = u\sqrt{p}\,\bar{r}\left(1 - \frac{\sqrt{p}\,\bar{r}}{m}\right) - \frac{kpc\bar{r}^2}{p(1+\bar{r}^2)} - qE\sqrt{p}\,\bar{r}. \quad (5)$$

Set $\beta = \frac{m}{\sqrt{p}} \Rightarrow m = \beta\sqrt{p},$

$$\frac{\sqrt{p}\,d\bar{r}}{\beta\sqrt{p}\,d\bar{t}} = \frac{u\bar{r}\sqrt{p}}{u\sqrt{p}}\left(1 - \frac{\bar{r}\sqrt{p}}{\beta\sqrt{p}}\right) - \frac{kc\bar{r}^2}{u\sqrt{p}(1+\bar{r}^2)} - \frac{qE\bar{r}\sqrt{p}}{u\sqrt{p}} \quad (6)$$

$$\frac{d\bar{r}}{d\bar{t}} = \beta\bar{r}\left(1 - \frac{\bar{r}}{\beta}\right) - \frac{\beta kc\bar{r}^2}{u\sqrt{p}(1+\bar{r}^2)} - \frac{\beta qE\bar{r}}{u} \quad (7)$$

$$\frac{d\bar{r}}{d\bar{t}} = \beta\bar{r} - \frac{\beta\bar{r}^2}{\beta} - \frac{kmc\bar{r}^2}{u\sqrt{p}^2(1+\bar{r}^2)} - \frac{mqE\bar{r}}{u\sqrt{p}} \quad (8)$$

$$\frac{d\bar{r}}{d\bar{t}} = \beta\bar{r} - \bar{r}^2 - \frac{kmc\bar{r}^2}{up(1+\bar{r}^2)} - \frac{mqE\bar{r}}{u\sqrt{p}} \quad (9)$$

$$\frac{d\bar{r}}{d\bar{t}} = \bar{r}(\beta - \bar{r}) - \frac{\alpha c\bar{r}}{(1+\bar{r}^2)} - \delta qE\bar{r} \quad (10)$$

where

$$\alpha = \frac{km}{u\sqrt{p}}, \qquad \delta = \frac{m}{u\sqrt{p}}$$

Removing bars, we have

$$\frac{dr}{dt} = r(\beta - r) - \frac{\alpha cr^2}{(1+r^2)} - \delta qEr \quad (11)$$

Similarly, from Equation (2) above, setting

$$\bar{c} = \frac{c}{\sqrt{p}} \Rightarrow c = \bar{c}\sqrt{p}, \quad dc = \sqrt{p}\,d\bar{c},$$

we have

$$\frac{d\bar{c}}{d\bar{t}} = v\bar{c}\left(1 - \frac{\bar{c}\sqrt{p}}{n}\right) + \frac{ekp\bar{r}^2\bar{c}}{p + p\bar{r}^2} - d\bar{c} \quad (12)$$

Setting

$$\bar{t} = \frac{vt\sqrt{p}}{n} \Rightarrow t = \frac{\bar{t}n}{v\sqrt{p}}, \qquad dt = \frac{n\,d\bar{t}}{v\sqrt{p}},$$

we have

$$\frac{d\bar{c}}{d\bar{t}} = \frac{n\bar{c}}{\sqrt{p}} - \bar{c}^2 + \frac{nek\bar{r}^2\bar{c}}{v\sqrt{p}(1+\bar{r}^2)} - \frac{nd\bar{c}}{v\sqrt{p}} \quad (13)$$

Set

$$\sigma = \frac{n}{\sqrt{p}} \Rightarrow n = \sigma\sqrt{p},$$

so that

$$\frac{d\bar{c}}{dt} = \bar{c}(\sigma - \bar{c}) + \frac{\rho \bar{r}^2 \bar{c}}{(1+\bar{r}^2)} - \mu \bar{c} \qquad (14)$$

Where

$$\rho = \frac{enk}{v\sqrt{p}}, \qquad \mu = \frac{dn}{v\sqrt{p}}$$

Removing bars, we have

$$\frac{dc}{dt} = c(\sigma - c) + \frac{\rho r^2 c}{(1+r^2)} - \mu c. \qquad (15)$$

Adding noise and periodic force terms, the model which is Equation (11) and (15) becomes

$$\frac{\partial r}{\partial t} = r(\beta - r) - \frac{\alpha r^2 c}{(1+r^2)} - \delta q E r + A\sin\omega t + D_1 \nabla^2 r \qquad (16)$$

$$\frac{\partial c}{\partial t} = c(\sigma - c) + \frac{\rho r^2 c}{(1+r^2)} - \mu c + \bar{A}\sin(\omega t + \varphi) + D_2 \nabla^2 c \qquad (17)$$

Where the non-negative constants $D_1$ and $D_2$ are respectively prey and predator diffusion coefficients. $\nabla^2 = \frac{\partial^2}{\partial x^2} + \frac{\partial^2}{\partial y^2}$ is the Laplacian operator in two–dimensional space.

The periodic force is assumed to be sinusoidal with amplitude A and frequency $\omega$ in equation (16). In equation (17), the periodic force and noise is also sinusoidal with amplitude denoted as $\bar{A}$, frequency $\omega$ and phase shift $\varphi$. This periodic force is considered to be positive reason being that the toxins produced by different populations have significant roles in shaping the dynamical behavior of ecosystems.

In this model, we called the noise in Equation (16) White noise because it has equal intensity at different frequencies and Colored noise is found in Equation (17) because it is closer to physical reality and have been used in describing ecological evolution. The noise in Equation (17) will have an influence over that of Equation (16).

*Existence of Equilibria*
From (16) and (17) above, that is, the Predator-Prey Model, set

$$\frac{\partial R}{\partial t} = 0, \frac{\partial C}{\partial t} = 0 \text{ and } D_1 = D_2 = 0$$

resulting in a system of non-linear algebraic equations to solve.

Putting $R_e$ and $C_e$ as the equilibrium solutions for the Rat and Cat populations respectively, the system of algebraic equations to solve from Equation (16) then becomes:

$$\beta R_e - R_e^2 - \frac{\alpha R_e^2 C_e}{1+R_e^2} - \delta q E R_e + A\sin\omega t = 0 \qquad (18)$$

Ignoring the H.O.T of $R_e$ and setting $C_e = 0$, Equation (18) gives

$$\beta R_e - \delta q E R_e + A\sin\omega t = 0. \qquad (19)$$

$$R_e = -\frac{A\sin\omega t}{(\beta - \delta q E)} \qquad (20)$$

Implying that we have possible equilibria at $C_e = 0$ or $R_e = -\frac{A\sin\omega t}{(\beta - \delta q E)}$

Similarly, from the Equation 17 above, we have;

$$C_e \sigma - C_e^2 + \frac{\rho R_e^2 C_e}{(1+R_e^2)} - \mu C_e + \bar{A}\sin(\omega t + \varphi) = 0 \qquad (21)$$

Ignoring the H.O.T. of $C_e$ and setting $R_e = 0$, Equation (21) gives

$$C_e(\sigma - \mu) = -\bar{A}\sin(\omega t + \varphi). \qquad (22)$$

$$C_e = -\frac{\bar{A}\sin(\omega t + \varphi)}{(\sigma - \mu)}. \qquad (23)$$

Implying that we have possible equilibria at $R_e = 0$ or $C_e = -\dfrac{\overline{A}\sin(\omega t + \varphi)}{(\sigma - \mu)}$

From the above, the simultaneous solution of the Equation (16) and (17) shows that when $R_e = 0$, then $C_e = 0$ which gives rise to the trivial solution (extinction of both species), that is,

$$(R_e, C_e) = (0,0).  \qquad (24)$$

The other equilibrium solution which denotes the co-existence of both species is given by:

$$(R_e, C_e) = \left(-\dfrac{A\sin\omega t}{(\beta - \delta q E)}, -\dfrac{\overline{A}\sin(\omega t + \varphi)}{(\sigma - \mu)}\right). \qquad (25)$$

Note here that, the equilibria do not help explain the oscillatory behavior of the data reflected by the Rat and Cat. So we need information about the stability of these equilibria before we can demonstrate that this is an appropriate model. Hence we perform a linear analysis.

**Linear/Stability analysis of the model**

In solving for the stability, we perform a linearization using partial derivatives. Thus, the Jacobian matrix of the Predator – Prey Model becomes:

$$J(R,C) = \begin{pmatrix} \dfrac{\partial U}{\partial R} & \dfrac{\partial U}{\partial C} \\ \dfrac{\partial V}{\partial R} & \dfrac{\partial V}{\partial C} \end{pmatrix}$$

Let $U = \beta R - R^2 - \dfrac{\alpha R^2 C}{1+R^2} - \delta q E R + A\sin\omega t$

$V = C\sigma - C^2 + \dfrac{\rho R^2 C}{1+R^2} - \mu C + \overline{A}\sin(\omega t + \varphi)$

$$J(R,C) = \begin{pmatrix} \beta - 2R - \dfrac{(1+R^2)(2\alpha RC) + 2\alpha R^3 C}{(1+R^2)(1+R^2)} - \delta q E & \dfrac{-\alpha R^2}{1+R^2} \\ \dfrac{(1+R^2)(2\rho RC) - 2\rho R^3 C}{(1+R^2)^2} & \sigma - 2C + \dfrac{\rho R^2}{1+R^2} - \mu \end{pmatrix} = \begin{pmatrix} R \\ C \end{pmatrix}$$

When evaluated at the steady state of (0, 0) the Jacobian matrix J is

$$J(0,0) = \begin{pmatrix} \beta - \delta q E & 0 \\ 0 & \sigma - \mu \end{pmatrix} - \lambda \begin{pmatrix} 1 & 0 \\ 0 & 1 \end{pmatrix} = \begin{pmatrix} \beta - \delta q E - \lambda & 0 \\ 0 & \sigma - \mu - \lambda \end{pmatrix}$$

The Eigen values are:

$\lambda_1 = \beta - \delta q E$ and $\lambda_2 = \sigma - \mu$

For $\lambda_1 > \lambda_2 > 0$

Hence, the system is unstable.

For the second solution of equilibrium point:

$$J\left(-\dfrac{A\sin\omega t}{\beta - \delta q E}, -\dfrac{\overline{A}\sin(\omega t + \varphi)}{\sigma - \mu}\right)$$

$$= \begin{pmatrix} \beta - 2R - \dfrac{2\alpha RC - 2\alpha R^3 C}{1+R^2} - \delta q E, & \dfrac{\alpha R^2}{1+R^2} \\ \dfrac{2\rho RC - 2\rho R^3 C}{1+R^2} & \sigma - 2C + \dfrac{\rho R^2}{1+R^2} - \mu \end{pmatrix} -$$

Put

$$\lambda \begin{pmatrix} 1 & 0 \\ 0 & 1 \end{pmatrix}$$

$A = \dfrac{R}{(1+R^2)}, \qquad B = \dfrac{R^2}{(1+R^2)}$

So that

$$J(R,C) = \begin{pmatrix} \beta - 2R - 2\alpha AC - 2\alpha BRC - \delta q E, & \alpha B \\ 2\rho AC - 2\rho BRC & \sigma - 2C + \rho B - \mu \end{pmatrix}$$

Where

$$P = \frac{A\sin\omega t}{\beta - \delta qE} = -R, \quad K = \frac{\overline{A}\sin(\omega t + \varphi)}{\sigma - \mu} = -C$$

$$J(R,C) = \begin{pmatrix} \beta + 2P + 2\alpha K(A - BP), & \alpha B \\ -2\rho K(A + BP), & \sigma + 2K + \rho B - \mu \end{pmatrix} - \begin{pmatrix} \lambda & 0 \\ 0 & \lambda \end{pmatrix}$$

The Eigen values are

$$\lambda_1 = (\beta + 2P + 2\alpha K(A - BP))(\sigma + 2K + \rho B - \mu)$$
and $\lambda_2 = -2\alpha BK(A + BP)$

For $\lambda_1 > 0, \lambda_2 < 0$;

Hence, the system is unstable.

**Simulation of the Model**

From Figures 1 to 2, the values of the parameters are: $\phi = -0.192$, $\beta = 0.2$, $\delta = 0.066$, $A = \overline{A} = 1$, $\omega = \frac{2\pi}{12}$, $Q = 2, M = 1, q = 1, E = 0.125, \sigma = 0.1$, $\rho = 0.05, \mu = 0.05, t = 0.05, r = 2 - 10$, $\theta = -0.05 - 0.05r^2$ (ignoring the H.O.T of $r > r^2$). $\varphi = \pm\frac{\pi}{4}$ (for increasing and decreasing phase shift).

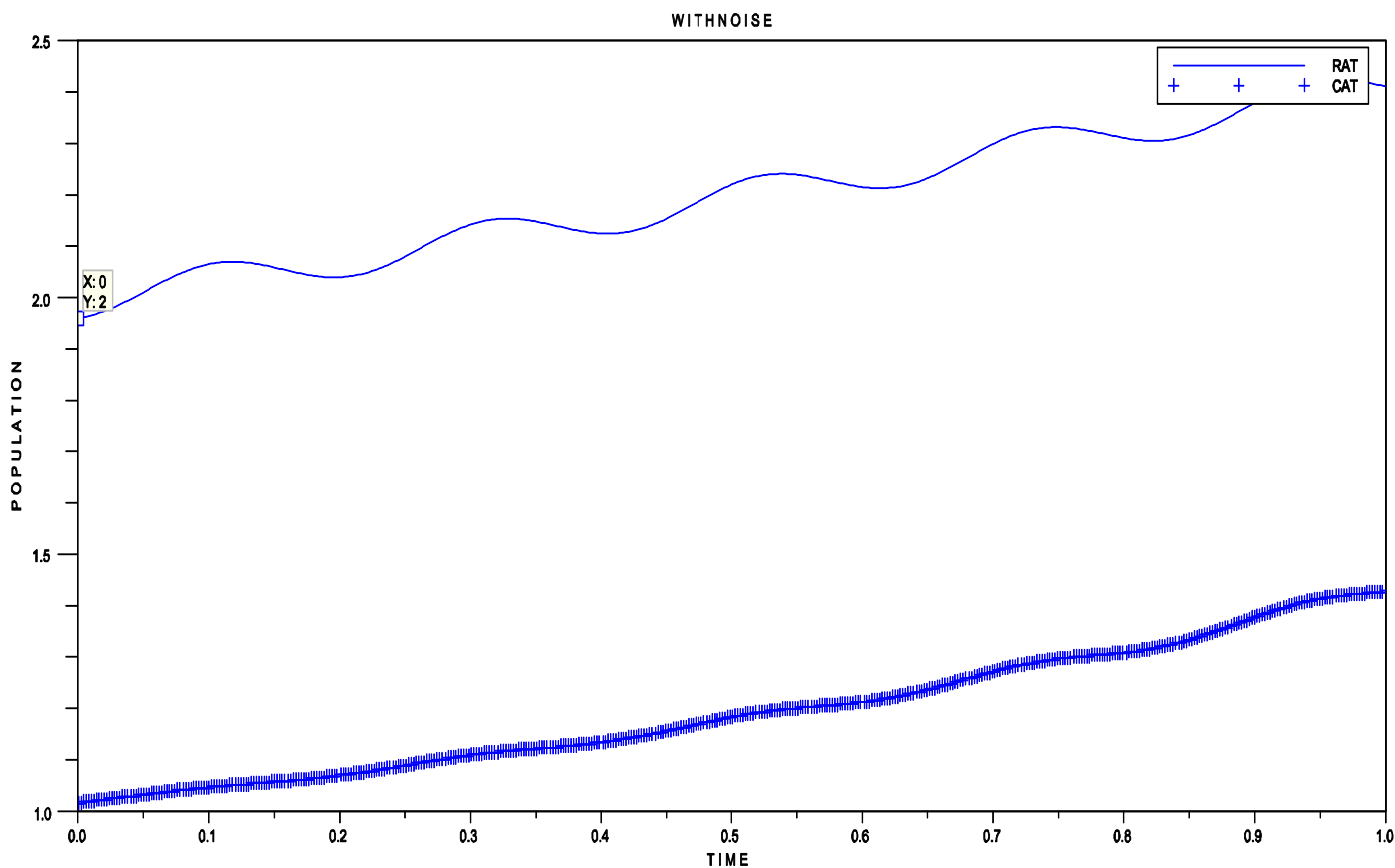

**Figure 1.** The Predator (Cat) and Prey (Rat) Population against Time with disturbances and increasing phase shift.

**DISCUSSION**

The outcome of the simulation showed that disturbances can cause oscillatory wave pattern. The model when initiated with different initial values lead to different curves of typical time series of predator C(t) and of prey R(t). As a result, solutions show oscillations with a frequency $\omega = \frac{2\pi}{12}$. We can say here that the plausible reason for the origin of the oscillations is delayed predator-prey interactions and the emergence of disturbance. Nevertheless, as the intensity of the disturbances are increased, the fluctuations become more apparent. It is understood that the frequency of the waves and resounding pattern arises when the disturbances and external periodic forces are present in the system (Equations 16 and 17). The results depict a regular frequency pattern which is a general characteristic of Predator-Prey Models. The figure shows the population of each species with time; but in contrast, both reveal that as the predator eats up the prey and its population density reduces at some time interval, the predator's population starts increasing at some time interval; also, when there is no much prey for the predator to eat its number starts declining while the prey's population density starts to climb and this happening is in a continuous process. The graph of the prey's population against time exhibits a higher, stronger and clearer sinusoidal wave pattern than that of the predator's population against time.

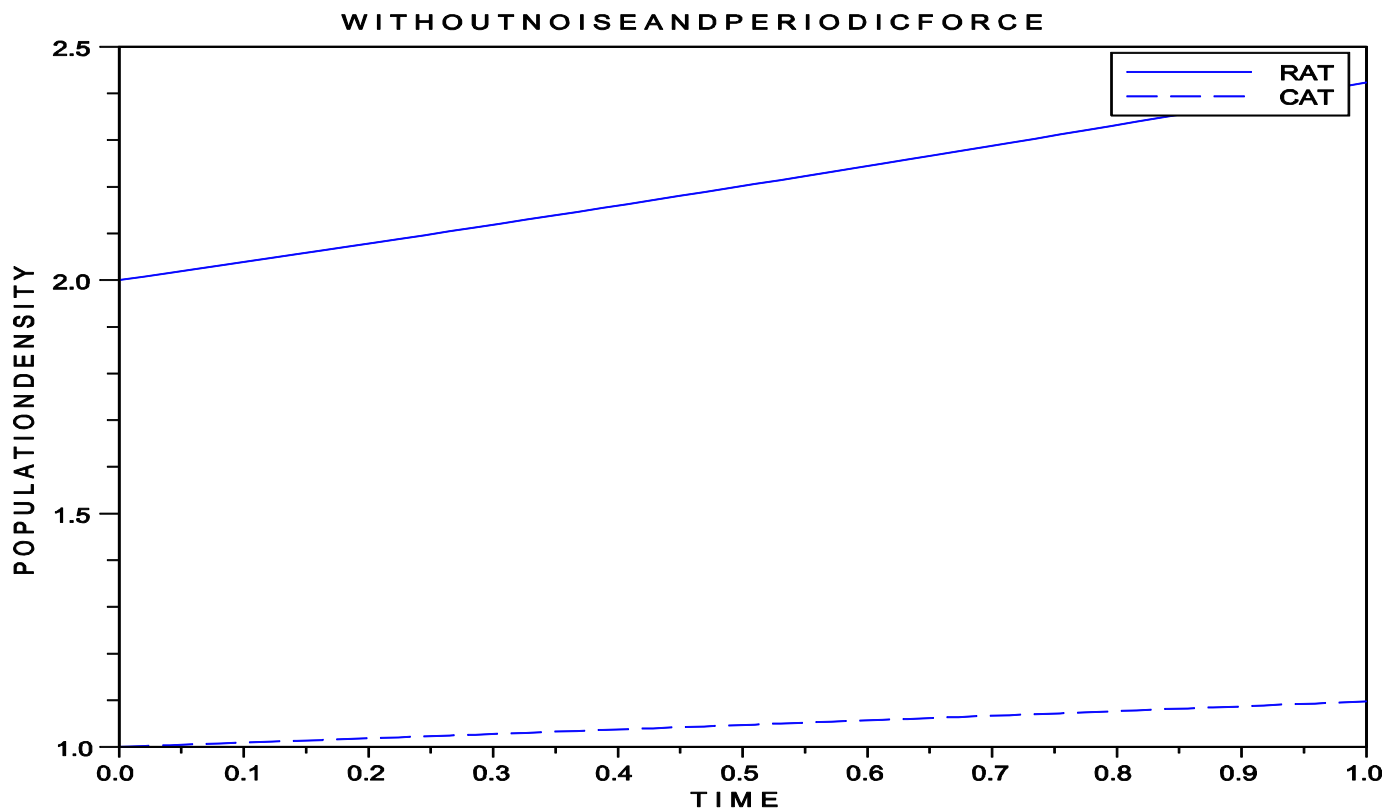

**Figure 2.** The Predator (Cat) and Prey (Rat) Population against Time with disturbances (Human) but without noise and periodic force.

The disturbances have an extensive effect in the model in that it supports the modulating nature of the species population especially that of the prey. Its presence in the system decreases the population density of the prey even in the absence of the predator or when the predator has not been introduced into the competition or system. Evidently, it is the presence of the disturbances that makes Figure 1 show a spiraling rise and fall nature.

**Conclusion**

This work has modified the Lotka-Volterra Predator-Prey Model by incorporating variables in the dynamics of the rat and cat population. Through simulation, this work has showed that disturbances in human form and noise when induced into the predator-prey system during competition, affects the interaction of both species. It can lead to the death of the prey causing a reduction in its population. The disturbances (especially noise) can also quarantine the prey in its hideout, making the predator to starve to death. The effect of the disturbances on the population of both species does not occur at the same time. At the time the effect is causing a reduction in the population of the prey specie, the predator population is increasing. Similarly, at the time the effect of the disturbances is causing a reduction in the population of the predator, the population of the prey is rising. The modified model is suitable for studying the dynamics of interspecific interaction between predators and their prey with disturbances.

**CONFLICT OF INTERESTS**

The authors have not declared any conflict of interests.